\documentstyle[12pt]{article}

\def\fr{\frac}

\def\de{\delta}

\def\a{a^{\dagger}}
\def\e{\exp_q}

\newcommand{\bd}{\begin{displaymath}}
\newcommand{\ed}{\end{displaymath}}
\newcommand{\bb}{\begin{equation}}
\newcommand{\ee}{\end{equation}}
\newcommand{\bea}{\begin{eqnarray}}
\newcommand{\eea}{\end{eqnarray}}

\begin{document}
\baselineskip 1.5 \baselineskip
%% Title Page 

\vspace{.2cm}

\begin{center}
\Large {\bf $gl_q(n)$-Covariant Multimode Oscillators and  q-Symmetric States 
}
\\[1cm]
\large W.-S.Chung \\[.3cm]
\normalsize  
Department of Physics and Research Institute of Natural Science,  \\
\normalsize  Gyeongsang National University,   \\
\normalsize   Jinju, 660-701, Korea
\end{center}

\vspace{0.5cm}
\begin{abstract}
In this paper the coherent states and q-symmetric states for
$gl_q(n)$-covariant multimode 
oscillator system are investigated. 
\end{abstract}

\setcounter{page}{1}
\section{Introduction}

Quantum                                                               
groups                                                                      
or                        
q-deformed                    
Lie         
algebra                
implies   
some    
specific 
deformations 
of                    classical                      Lie           algebras.

From           
a             
mathematical      
point   
of  
view,     
it 
is  
a 
non-commutative                                                  
associative                                                                 
Hopf         
algebra.         
The       
structure        
and 
representation 
theory    of                                   quantum               groups  
have                 
been             
developed        
extensively  
by  
Jimbo    
[1] 
and 
Drinfeld                                                                
[2].                                                                        
            
The                                   
q-deformation                    
of        
Heisenberg       
algebra   
was     
made  
by 
Arik and Coon [3], Macfarlane [4] and Biedenharn [5].
Recently                                                               
there                                               
has                           
been                    
some              
interest               
in    
more     
general   
deformations 
involving                                                                 
an                                                                          
arbitrary                      
real                           
functions     
of           
weight    
generators  
and   
including 
q-deformed algebras as a special case [6-10].

\def\a{a^{\dagger}}

\def\q{q^{-1}}
\def\a{a^{\dagger}}
Recently  Greenberg [11]  has  studied the  following  q-deformation of
multi  mode boson 
algebra:
\bd
a_i \a_j -q \a_j a_i=\de_{ij},
\ed
where the deformation parameter $q$ has to be real.
The main problem of Greenberg's approach is that we can not derive the relation
among $a_i$'s operators at all.
Moreover the above algebra is not  covariant under $gl_q(n)$ algebra.
In order to solve this problem we should find the q-deformed multimode
oscillator
algebra which is covariant under $gl_q(n)$ algebra.
Recently the Fock space representation of  $gl_q(n)$-covarinat multimode
oscillator
system was known by some authers [12].
In this paper we construct the correct form of coherent states for the above 
mentioned oscillator  system and  obtain the  q-symmetric states
generalizing  the bosonic 
states.

\section{Coherent states of $gl_q(n)$-covariant oscillator algebra}

$gl_q(n)$-covariant oscillator algebra is defined as [12]
\bd
\a_i\a_j =q \a_j \a_i,~~~(i<j)  
\ed
\bd
a_ia_j=\fr{1}{q}a_j a_i,~~~(i<j)
\ed
\bd
a_i\a_j=q \a_ja_i,~~~(i \neq j)
\ed
\bd
a_i\a_i =1+q^2 \a_ia_i +(q^2-1) \sum_{k=i+1}^n\a_k a_k,~~~(i=1,2,\cdots,n-1)
\ed
\bd
a_n \a_n =1+q^2 \a_n a_n,
\ed
\bb
[N_i, a_j]=-\de_{ij}a_j,~~~[N_i, \a_j]=\de_{ij}\a_j,~~~(i,j=1,2,\cdots, n )
\ee
where we restrict our concern to the case that $q$ is real and $0<q<1$.
Here $N_i$ plays a role of number operator and $a_i(\a_i)$ plays a role of
annihilation(creation) operator.
From the above algebra one can obtain the relation between the number operators
and mode opeartors as follows
\bb
\a_ia_i=q^{2\sum_{k=i+1}^nN_k}[N_i],
\ee
where $[x]$ is called a q-number and is defined as
\bd
[x]=\fr{q^{2x}-1}{q^2-1}.
\ed
\def\nn{|n_1,n_2,\cdots,n_n>}
Let us introduce the Fock space basis $\nn$ for the number operators
$N_1,N_2,\cdots, N_n$ 
satisfying
\bb
N_i\nn=n_i\nn,~~~(n_1,n_2,\cdots,n_n=0,1,2\cdots)
\ee
Then we have the following representation
\bd
a_i\nn=q^{\sum_{k=i+1}^nn_k}\sqrt{[n_i]}|n_1,\cdots, n_i-1,\cdots,n_n>
\ed
\bb
\a_i\nn=q^{\sum_{k=i+1}^nn_k}\sqrt{[n_i+1]}|n_1,\cdots, n_i+1,\cdots,n_n>.
\ee

From the above representation we  know that there exists the ground state
$|0,0,\cdots,0>$ 
satisfying
$a_i|0,0>=0$ for  all $i=1,2,\cdots,n$. Thus the  state $\nn$ is obtatind
by applying the 
creation operators
to the ground state $|0,0,\cdots,0>$
\bb
\nn=\fr{(\a_n)^{n_n}\cdots(\a_1)^{n_1}}{\sqrt{[n_1]!\cdots
[n_n]!}}|0,0,\cdots,0>.
\ee
If we introduce the scale operators as follows
\bb
Q_i=q^{2N_i},~~(i=1,2,\cdots,n),
\ee
we have from the algebra (1)
\bb
[a_i,\a_i]=Q_iQ_{i+1}\cdots Q_n.
\ee
Acting the operators $Q_i$'s on the basis $\nn$ produces
\bb
Q_i\nn=q^{2n_i}\nn .
\ee

From the relation  $a_i a_j =\fr{1}{q}a_j a_i,~~(i<j)$, the  coherent states
for $gl_q(n)$ 
algebra
is defined as
\bb
a_i|z_1,\cdots,z_i,\cdots,z_n>=z_i|z_1,\cdots, z_{i},qz_{i+1},\cdots,q z_n>.
\ee
Solving the eq.(9) we obtain
\bb
|z_1,z_2,\cdots,z_n>=c(z_1,\cdots,z_n)\sum_{n_1,n_2,\cdots,n_n=0}^{\infty} 
\fr{z_1^{n_1}z_2^{n_2}\cdots z_n^{n_n}}{\sqrt{[n_1]![n_2]!\cdots [n_n]!}}\nn .
\ee
Using eq.(5) we can rewrite eq.(10) as
\bb
|z_1,z_2,\cdots,z_n>=c(z_1,\cdots,z_n) 
\e(z_n\a_n)\cdots\e(z_2\a_2)\e(z_1\a_1)|0,0,\cdots,0>.
\ee
where q-exponential function is defined as
\bd
\e(x)=\sum_{n=0}^{\infty}\fr{x^n}{[n]!}.
\ed
The q-exponential function satisfies the following recurrence relation
\bb
\e(q^2 x)=[1-(1-q^2)x]\e(x)
\ee
Using the above relation and the fact that $0<q<1$, we obtain the formula
\bb
\e(x) =\Pi_{n=0}^{\infty}\fr{1}{1-(1-q^2)q^{2n}x}
\ee
Using the normalization of the coherent state , we have
\bb
c(z_1,z_2,\cdots,z_n)=\e(|z_1|^2)\e(|z_2|^2)\cdots \e(|z_n|^2).
\ee
The coherent state satisfies the completeness relation
\bb
\int\cdots  \int
|z_1,z_2,\cdots,z_n><z_1,z_2,\cdots,z_n|\mu(z_1,z_2,\cdots,z_n)  d^2z_1 
d^2z_2\cdots d^2 z_n=I,
\ee
where the weighting function $\mu(z_1,z_2,\cdots,z_n)$ is defined as
\bb
\mu(z_1,z_2,\cdots,z_n)=\fr{1}{\pi^2}\Pi_{i=1}^n\fr{\e(|z_i|^2)}
{\e(q|z_i|^2)}.
\ee
In deriving eq.(15) we used the formula
\bb
\int_0^{1/(1-q^2)}x^n \e(q^2 x)^{-1} d_{q^2} x=[n]!
\ee

\def\ot{\otimes}
\section{q-symmetric states}

In this section we study the statistics of many particle state.
Let $N$ be the number of particles. Then the N-partcle state can be obtained
from
the tensor product of single particle state:
\bb
|i_1,\cdots,i_N>=|i_1>\ot |i_2>\ot \cdots \ot |i_N>,
\ee
where $i_1,\cdots, i_N$ take one value among $\{ 1,2,\cdots,n \}$ and the sigle
particle state is defined by $|i_k>=\a_{i_k}|0>$.

Consider the case that k appears $n_k$ times in the set $\{ i_1,\cdots,i_N\}$.
Then we have
\bb
n_1 + n_2 +\cdots + n_n =\sum_{k=1}^n n_k =N.
\ee
Using these facts we can define the q-symmetric states as follows:
\bb 
|i_1,\cdots, i_N>_q
=\sqrt{\frac{[n_1]!\cdots [n_n]!}{[N]!}}
\sum_{\sigma \in Perm}
\mbox{sgn}_q(\sigma)|i_{\sigma(1)}\cdots i_{\sigma(N)}>,
\ee
where
\begin{displaymath}
\mbox{sgn}_q(\sigma)=
q^{R(i_1\cdots i_N)}q^{R(\sigma(1)\cdots \sigma(N))},
\end{displaymath}
\bd
R(i_1,\cdots,i_N)=\sum_{k=1}^N\sum_{l=k+1}^N R(i_k,i_l)
\ed
and
\bd
R(i,j)=\cases{
1 & if $ i>j$ \cr
0 & if $ i \leq j $ \cr
}
\ed
Then the q-symmetric states obeys
\begin{equation}
|\cdots, i_k,i_{k+1},\cdots>_q=
\cases{
q^{-1} |\cdots,i_{k+1},i_k,\cdots>_q& if $i_k<i_{k+1}$\cr
 |\cdots,i_{k+1},i_k,\cdots>_q& if $i_k=i_{k+1}$\cr
q |\cdots,i_{k+1},i_k,\cdots>_q& if $i_k>i_{k+1}$\cr
}
\end{equation}
The above property can be rewritten by introducing the deformed transition
operator
$P_{k,k+1}$ obeying
\bb
P_{k,k+1}
|\cdots, i_k , i_{k+1},\cdots>_q =|\cdots, i_{k+1},i_k,\cdots>_q
\ee
This operator satisfies 
\bb
P_{k+1,k}P_{k,k+1}=Id,~~~\mbox{so}~~P_{k+1,k}=P^{-1}_{k,k+1}
\ee
Then the equation (21) can be written as
\bb
P_{k,k+1}
|\cdots, i_k , i_{k+1},\cdots>_q 
=q^{-\epsilon(i_k,i_{k+1})}
|\cdots, i_{k+1},i_k,\cdots>_q
\ee
where $\epsilon(i,j)$ is defined as
\bd
\epsilon(i,j)=
\cases{
1 & if $ i>j$\cr
0 & if $ i=j$ \cr
-1 & if $ i<j$ \cr }
\ed
The relation (24) goes to the symmetric relation for the ordinary bosons 
 when the deformation parameter $q$ goes to $1$.
If we define the fundamental q-symmetric state $|q>$ as
\bd
|q>=|i_1,i_2,\cdots,i_N>_q
\ed
with $i_1 \leq i_2 \leq \cdots \leq i_N$, we have 
\bd
||q>|^2 =1.
\ed
In deriving the above relation we used following identity
\bb
\sum_{\sigma \in Perm } q^{2R(\sigma(1),\cdots, \sigma(N))}=
\fr{[N]!}{[n_1]!\cdots [n_n]!}
\ee
The derivation of above formula will be given in Appendix.

\section{Concluding Remark}
In this paper  the $gl_q(n)$-covariant oscillator algebra and its
coherent states are discussed. The q-symmetric states generalizing the
symmetric 
(bosonic) states are obtained by using the $gl_q(n)$-covariant oscillators
and are shown to be orthonormal. 
I think that the q-symmetric states will be important when we consider the new
statistical field theory generalizing the ordinary one.

\section*{Appendix}
In this appendix we prove the relation(25) by using the mathematical induction.
Let us assume that the relation (25) holds for $N$.
Now we should prove that eq.(25) still hold for $N+1$.
Let us consider the case that $i$ appears $n_i+1$ times.
Then we should show
\def\s{\sigma}
\bb
\sum_{\s \in Perm}q^{2R(\s(1), \cdots, \s(N+1))}=\fr{[N+1]!}{[n_1]!\cdots 
[n_{i-1}]![n_i+1]![n_{i+1}]!\cdots [n_n]!}
\ee
In this case the above sum can be written by three pieces:
\bb
\sum_{j=1}^{i-1} \sum_{ \s(1)=j}
+\sum_{\s(1)=i}
+\sum_{j=i+1}^{n} \sum_{ \s(1)=j}
\ee
Thus the left hand side of eq.(26) is given by
\bea
LHS&=& \sum_{j=1}^{i-1}\sum_{\s(1)=j} q^{2R(j,\s(2),\cdots,\s(n+1))}\cr
&+& \sum_{\s(1)=i} q^{2R(i,\s(2),\cdots,\s(n+1))}                   \cr
&+& \sum_{j=i+1}^{n}\sum_{\s(1)=j} q^{2R(j,\s(2),\cdots,\s(n+1))}      \cr
\eea
Then we have
\bd
R(j,\s(2),\cdots,\s(n+1)) =\sum_{k=2}^{N+1} R(j,\s(k))
+R(\s(2),\cdots,\s(N+1))
\ed
where
\bd
\sum_{k=2}^{N+1}R(j,\s(k))
=\cases{ n_1 + \cdots + n_{j-1} & if $ j \leq i$\cr
n_1 +\cdots + n_{j-1} +1 & if $ j >i$ \cr}
\ed
Using the above relations the LHS of eq.(26) can be written as
\bea
LHS
&=& \sum_{j=1}^{i-1}q^{2(n_1  + \cdots  + n_{j-1})}
\fr{[N]!}{[n_i+1]![n_j-1]!\Pi_{k \neq 
i,j} [n_k]!}    \cr
&+& q^{2(n_1 + \cdots + n_{i-1})} \fr{[N]!}{\Pi_{k } [n_k]!}    \cr
&+&  \sum_{j=i+1}^{N+1}q^{2(n_1 +  \cdots  +  n_{j-1}+1)}
\fr{[N]!}{[n_i+1]![n_j-1]!\Pi_{k 
\neq i,j} [n_k]!}    \cr
\eea
If we pick up the common factor of three terms of eq.(29), we have
\bd
I=J
\fr{[N]!}{[n_i+1]!\Pi_{k \neq i} [n_k]!}  
\ed
where
\bea
J&=&\fr{1}{q^2-1}[
\sum_{j=1}^{i-1}q^{2(n_1 + \cdots + n_{j-1})} [n_j]
+q^{2(n_1 + \cdots + n_{i-1})}[n_1+1] 
+\sum_{j=i+1}^{N+1}q^{2(n_1 + \cdots + n_{j-1}+1)}[n_j] \cr
&=&[N+1]
\eea
Thus we proved the relation (25).

\section*{Acknowledgement}
This                   paper                was
supported         by  
the   KOSEF (961-0201-004-2)   
and   the   present   studies    were   supported   by   Basic  
Science 
Research Program, Ministry of Education, 1995 (BSRI-95-2413).

\vfill\eject

\end{document}